\RequirePackage{fix-cm}
\documentclass[smallextended]{svjour3}       
\smartqed  
\usepackage{graphicx} 
\usepackage{extarrows}
\usepackage{latexsym}
\usepackage{amsmath}
\usepackage{amssymb}
\usepackage{amsfonts}
\usepackage{verbatim}
\usepackage{mathrsfs}
\usepackage[modulo]{lineno} 
\usepackage[latin1]{inputenc}
\usepackage{color}
\usepackage{xcolor}
\usepackage[colorlinks,citecolor=blue,urlcolor=blue]{hyperref}
\usepackage{soul}
\usepackage{enumitem}

\newtheorem{assumption}{Assumption}

\numberwithin{equation}{section}

\newcommand{\ee}{\mathbb E}

\newcommand{\pp}{\mathbb P}
\newcommand{\nn}{\mathbb N}
\newcommand{\rr}{\mathbb R}

\newcommand{\BB}{\mathcal B}

\newcommand{\FFF}{\mathscr F}

\newcommand{\<}{\langle}
\renewcommand{\>}{\rangle}
\allowdisplaybreaks \allowdisplaybreaks[4]

\newcommand{\dd}{\mathrm{d}}

\newcommand{\abs}[1]{\left\lvert #1 \right\rvert}
\newcommand{\norm}[1]{\left\lVert #1 \right\rVert}

\begin{document}

\title{Geometric Ergodicity and Strong Error Estimates for Tamed Schemes of Super-linear SODEs 
}
\subtitle{}

\titlerunning{Geometric Ergodicity and Error Estimates of Tamed Schemes Super-linear SODEs}        

\author{Zhihui LIU         \and
        Xiaoming WU 
}


\institute{Z. LIU \at
              Department of Mathematics \& National Center for Applied Mathematics Shenzhen (NCAMS) \& Shenzhen International Center for Mathematics, Southern University of Science and Technology, Shenzhen, 518055, P.R.China  \\
              \email{liuzh3@sustech.edu.cn (Corresponding author)}           
           \and
           X. WU \at
              Department of Mathematics, Southern University of Science and Technology, Shenzhen 518055, P.R. China \\
              \email{12331004@mail.sustech.edu.cn}
}

\date{Received: date / Accepted: date}

\maketitle

\begin{abstract}
We construct a family of explicit tamed Euler--Maruyama (TEM) schemes, which can preserve the same Lyapunov  function for super-linear stochastic ordinary differential equations (SODEs) driven by multiplicative noise.
These TEM schemes are shown to inherit the geometric ergodicity of the considered SODEs and converge with optimal strong convergence orders. 
Numerical experiments verify our theoretical results.

\end{abstract}


\section{Introduction}

The long-time behavior of the Wiener process-driven SODE 
\begin{align}\label{sde}
    \dd X(t) = b(X(t))\dd t + \sigma(X(t)) \dd W(t), ~ t \ge 0, \tag{SDE}
\end{align}
  plays a vital role in many scientific areas.
As a significant long-time behavior, the ergodicity characterizes the identity of the temporal average and spatial average for a Markov process or chain generated by Eq. \eqref{sde} and its numerical discretization, respectively, which has a lot of applications in quantum mechanics, fluid dynamics, financial mathematics, and many other fields \cite{DZ96,HW19}. 

It is known that the coefficients of most nonlinear SODEs arising in applications do not satisfy the traditional, but restrictive, Lipschitz continuity assumptions.  
This paper analyzes integrators for the super-linear Eq. \eqref{sde} whose solution is uniquely ergodic with respect to an equilibrium distribution  that
\begin{enumerate}
\item
are ergodic with respect to an invariant measure on infinite time intervals, 

\item
strongly converge with optimal convergence rate to the solutions of Eq. \eqref{sde} on any finite time intervals, and

\item

involve negligible computational expense. 
\end{enumerate} 
It is known, even for the particular Langevin system, that pure sampling methods can accomplish (1),  but they typically do not approximate the solution to Eq. \eqref{sde}; see, e.g., \cite{BV10}.
Integrators for Eq. \eqref{sde} certainly satisfy (2), but they are often divergent on infinite-time intervals or ergodic concerning a different equilibrium distribution \cite{MSH02,MT21}.
The backward Euler method and the stochastic theta method (see, e.g., \cite{LMW23,Liu24,LL25,MSH02}) satisfy (1) and (2), but they usually entail high computational cost.
We will show that a family of TEM methods constructed in this paper, as explicit schemes, can simultaneously accomplish these three goals.
Similar arguments were successfully applied in \cite{LS25} to construct Galerkin-based fully discrete tamed schemes for superlinear stochastic PDEs.

To construct explicit methods that inherit the unique ergodicity of Eq. \eqref{sde}, including the tamed and truncated methods studied, e.g., in \cite{BDMS19,HJK12,LNSZ24,NNZ25}, the numerical Lyapunov structure plays a key role.
However, it was shown in \cite{HJK11} that the classical Euler--Maruyama (EM) scheme applied to Eq. \eqref{sde} with super-linear growth coefficients would blow up in $p$-th moment for all $p \ge 2$. 
In particular, the square function is not an appropriate Lyapunov function of the EM scheme, even though it is a natural Lyapunov function of the considered monotone Eq. \eqref{sde}.

Our main aim in this paper is to construct a family of explicit TEM schemes to preserve the same Lyapunov structure of the super-linear Eq. \eqref{sde}.

 It is shown that, under certain non-degeneracy conditions ensuring the equivalence of transition probabilities, these TEM schemes inherit the geometric ergodicity of Eq. \eqref{sde}. 
To the best of our knowledge, this is the first time that explicit schemes have been constructed which preserve the unique ergodicity of super-linear SODEs driven by multiplicative noise. 
 To establish the second property, it is essential to estimate the effect on the dynamics of the TEM methods through strong error estimates between the TEM methods and  Eq. \eqref{sde}. 
 By examining the relationship between these error estimates and the step size, we prove that the expected optimal strong convergence order is $1/2$ in the multiplicative noise case and $1$ in the additive noise case.
 
The outline of the paper is as follows.
In Section \ref{sec2}, we propose the TEM schemes and develop their probabilistic regularity properties.  
The strong error estimates between the TEM schemes and the SODEs are explored in Section \ref{sec3}.
In Section \ref{sec4}, we validate the theoretical results by numerical experiments.

  \section{Geometric Ergodicity of Tamed Euler--Maruyama Schemes}
 \label{sec2}

This section will investigate the numerical ergodicity of explicit schemes applied to the $d$-dimensional Eq. \eqref{sde} driven by an $\rr^m$-valued Wiener process $W$ on a complete filtered probability space $(\Omega, \FFF, \mathbb{F}:=(\FFF(t))_{t\ge 0}, \pp)$ through exploring their Lyapunov structure and probabilistic regularity properties.
   
For the study of the Lyapunov structure of numerical schemes, 
our central assumption on the coefficients $b$ and $\sigma$ in \eqref{sde} is the following coupled coercivity and polynomial growth conditions.
Throughout, we denote by $\abs{\cdot}$ and $\<\cdot, \cdot\>$ the Euclidean norm and inner product in $\rr$, $\rr^d$, or $\rr^m$ if there is no confusion, and by $\norm{\cdot}$ the Hilbert--Schmidt norm in $\rr^{d \times m}$.

\begin{assumption} \label{ap}   
There exist positive constants $L_i$,  $i=0,1,2,3,4$, and $q$ such that for all $x \in \rr^d$, 
\begin{align} 
2 \<x-y,b(x)-b(y)\> + \|\sigma(x)-\sigma(y)\|^2 & \le L_0 |x-y|^{2}, \label{coe-coup} \\
2 \<x,b(x)\> + \|\sigma(x)\|^2 & \le L_1 - L_2 |x|^{q+2}, \label{coe} \\
|b(x)| & \le L_3+ L_4 |x|^{q+1}. \label{b-grow}
\end{align}  
 \end{assumption}

 Under Assumption \ref{ap} and certain integrability condition $b$ and $\sigma$, one can show the existence and uniqueness of the $(\FFF_t)_{t\ge 0}$-adapted solution to Eq. \eqref{sde}; see, e.g., \cite[Theorem 3.1.1]{LR15}.

\begin{remark} \label{rk-ap}
 \begin{enumerate}  
 \item
The coupled coercivity condition \eqref{coe} and the polynomial growth condition \eqref{b-grow} on $b$ imply the polynomial growth condition on $\sigma$, i.e., there exist positive constants $L_5$ and $L_6$ such that for all $x \in \rr^d$, 
 \begin{align}\label{s-grow+}
\|\sigma(x)\| ^2 \le L_5+L_6 |x|^{q+2}. 
 \end{align} 
 
\item 
 In particular, Assumption \ref{ap} encompasses polynomial drift (with a negative leading coefficient) and diffusion functions. 
 Indeed, let $k \in \nn_+$ and assume that $b:\rr^d \rightarrow \rr^d$ and $\sigma: \rr^d \rightarrow \rr^{d\times m}$ are polynomial functions, i.e., each component of $b$ and $\sigma$ is a polynomial in the coordinates $x_1,\cdots, x_d$. In the one-dimensional case $d=1$, this reduces to the familiar form
 \begin{align}\label{bs-ex}
b(x)=\sum_{i=0}^{2k+1} a_i x^i, \quad
\sigma(x)=\sum_{j=0}^{k+1} c_j x^j,
\quad x \in \rr,
 \end{align} 
 where $a_i, c_j \in \rr$, $i=0,1,\cdots,2k+1$, $j=0,1,\cdots,k+1$.  
Through elementary algebraic calculations, one derives \eqref{coe} and \eqref{b-grow} for positive constants $L_i$, $i=1,2,3,4$, depending on the coefficients of the polynomials in \eqref{bs-ex} and $q=2k$ provided 
 \begin{align}\label{bs-ex-con}
2a_{2k+1}+c_{k+1}^2<0.
 \end{align} 
In particular, this includes the double-well potential drift corresponding to $b(x) = (1-|x|^2) x$, $x \in \rr^d$, in the additive noise case ($\sigma(x) \equiv {\rm B} \in \rr^{d \times m}$, $x \in \rr^d$); see \cite[Remark 2.3 and Proposition 3.1]{NNZ25}.

\item

If $\sigma$ has at most linear growth (including the additive noise case),
 then the condition \eqref{coe} can be decomposed into
\begin{align} 
2 \<x,b(x)\> & \le L_1'- L_2' |x|^{q+2}, \label{b-coe} \\ 
\|\sigma(x)\|^2 & \le L_3'+ L_4' |x|^2, \label{s-grow}
\end{align}  
for all $x \in \rr^d$ with some positive constants $L_1'$, $L_2'$, $L_3'$, and $L_4'$, which had been studied in \cite{NNZ25} for the additive noise case.  
\end{enumerate}
\end{remark}

To construct explicit schemes of Eq. \eqref{sde}, as the coercive condition \eqref{coe} includes super-linear drift and diffusion coefficients, we utilize a taming technique that tames both the drift and diffusion coefficients.
More precisely, we consider the following tamed Euler--Maruyama (TEM) scheme applied to Eq. \eqref{sde}:
 \begin{align} \label{tem}
 	Z_n=Z_{n-1}+
 	b_\tau (Z_{n-1}) \tau 
 	+ \sigma_\tau(Z_{n-1})\delta_{n-1} W,
~ n \in \nn_+; \quad Z_0=X_0,
 \end{align}
 where $\tau \in (0, 1]$ denotes the temporal stepsize, $\delta_{n-1} W:=W(t_n) -W(t_{n-1})$, $n \in \nn_+$, and the tamed drift and diffusion functions are chosen such that $b_\tau: \rr^d \to \rr^d$ and $\sigma_\tau: \rr^d \to \rr^{d \times m}$   both have at most linear growth, and there exist positive constants $T_1$ and $T_2$,   independent of $\tau$, such that for all $x \in \rr^d$ and $\tau \in (0,1]$, 
\begin{align}  \label{coe-tau+}
2 \<x,b_\tau(x)\> + \|\sigma_\tau(x)\|^2 + \tau |b_\tau(x)|^2 
& \le T_1 - T_2 |x|^2. 
\end{align}  

We note that the continuous-time coercive condition \eqref{coe} inspires the condition \eqref{coe-tau+} but is naturally adapted to our discrete setting. 
Moreover, we assume that $b_\tau$ and $\sigma_\tau$ converge to $b$ and $\sigma$ when $\tau \rightarrow 0$ to make it consistent.
As the TEM scheme \eqref{tem} is explicit, it can be solved pathwise uniquely.
 Moreover, it is clear that $\{Z_n, \FFF_{t_n}\}_{n \in \nn}$ is a (time-homogeneous) Markov chain.

\subsection{Lyapunov Structure}

In this part, we investigate the Lyapunov Structure of the TEM scheme \eqref{tem}.

We begin with the following result about two examples such that $b_\tau$ and $\sigma_\tau$   both have at most linear growth and satisfy \eqref{coe-tau+}.
 
  \begin{lemma}
 Under Assumption \ref{ap}, there exists $\tau_{\max} \in (0, 1]$ such that for any $\tau \in (0, \tau_{\max}]$,  
 \begin{align} \label{bs-tau}
 b_\tau (x):=\frac{b(x)}{(1+\tau|x|^{2q})^{1/2}}, \quad 
 \sigma_\tau(x):=\frac{\sigma(x)}{(1+\sqrt{\tau}|x|^q)^{1/2}}, \quad 
x \in \rr^d,
\end{align}  
or 
\begin{align} \label{bs-tau+}
 b_\tau (x):=\frac{b(x)}{(1+\tau|x|^{2q})^{1/2}}, \quad 
 \sigma_\tau(x):=\frac{\sigma(x)}{(1+\sqrt{\tau}|x|^{2q})^{1/2}}, \quad 
x \in \rr^d,
\end{align}  
 both have at most linear growth and satisfy \eqref{coe-tau+}.
 \end{lemma}

 \begin{proof}
 It is clear that $b_\tau$ and $\sigma_\tau$ defined in \eqref{bs-tau} or \eqref{bs-tau+} have both at most linear growth.
 To show that \eqref{coe-tau+} is valid for \eqref{bs-tau} and \eqref{bs-tau+}, following \eqref{coe}, \eqref{b-grow}, and the elementary inequality $\sqrt{|a|+|b|}\le\sqrt{|a|}+\sqrt{|b|}$ and $(a+b)^2 \le 2a^2+2b^2 $ for all $a, b \in \rr$, we have 
 \begin{align*} 
& 2 \<x, b_\tau(x)\> + \|\sigma_\tau(x)\|^2 + \tau |b_\tau(x)|^2  \\
& \le \frac{2\<x, b(x)\>+\|\sigma(x)\|^2}{(1+\tau|x|^{2q})^{1/2}}
 + \frac{\tau | b(x)|^2}{1+\tau|x|^{2q}} \\
& =L_1+ 2 L_3^2 \tau
+ \Big(\frac{-L_2 |x|^q}{(1+\tau |x|^{2q})^{1/2}} 
+\frac{2L_4^2 \tau |x|^{2 q}}{1+\tau |x|^{2q}} \Big)|x|^2,
\quad x \in \rr^d.
 \end{align*} 
Let $\xi \ge 0$.
Then
 \begin{align} \label{est-con} 
& \frac{-L_2 \xi^q}{(1+\tau \xi^{2q})^{1/2}} 
+\frac{2L_4^2 \tau \xi^{2 q}}{1+\tau \xi^{2q}} \nonumber \\ 
& = - \frac{\beta \xi^q}{(1+\tau \xi^{2q})^{1/2}}
- \Big( \frac{L_2-\beta}{(1+\tau \xi^{2q})^{1/2}}
-\frac{2\tau L_4^2 \xi^q}
{1+\tau \xi^{2q}} \Big) \xi^q  \nonumber \\
&=: -\beta D_1^\tau(\xi^q)  - \xi^q D_2^\tau(\xi^q),
 \end{align}
 where $\beta$ is any constant in $(0, L_2)$ and  
 \begin{align} \label{d12}
 D_1^\tau(\eta) :=\frac{\eta}{(1+\tau \eta^2)^{1/2}},\quad
 D_2^\tau(\eta) :=\frac{L_2-\beta}{(1+\tau \eta^2)^{1/2}}
 - \frac{2\tau L_4^2 \eta}{1+\tau \eta^2},
\quad \eta \in \rr_+.
 \end{align} 
 
For any $\alpha>0$, it follows from the fact $D_1^\tau$ is non-decreasing for all $\tau \in (0, 1)$ that $(1 \ge )D_1^\tau(\eta) \ge \alpha/\sqrt{1+\alpha^2}$ for all $|\eta| \ge \alpha$.
This shows that 
\begin{align} \label{d1}
D_1^\tau(\xi^q) \ge \alpha^q/\sqrt{1+
\alpha^{2q}}, \quad \xi  > \alpha>0.
\end{align} 
On the other hand, applying the elementary inequality 
\begin{align} \label{in-sqrt}
\sqrt{x+y}\ge (\sqrt{x}+\sqrt{y})/\sqrt{2}, \quad \forall~ x, y \in \rr_+,
 \end{align}
we get
 \begin{align*}
 D_2^\tau(\eta) &=\frac{(L_2-\beta)(1+\tau\eta^2)^{1/2} -2 \tau L_4^2 \eta}{1+\tau \eta^2} 
 \ge \frac{\frac{L_2-\beta}{\sqrt{2}}(1+\sqrt{\tau}\eta) -2 \tau L_4^2 \eta}{1+\tau \eta^2}   \ge 0,  
 \end{align*}
 for all $\eta\in \rr_+$, where we have used the fact that $\sqrt \tau (L_2-\beta)/\sqrt{2} \ge 2 \tau L_4^2$ for all $\tau \in (0, \frac{(L_2-\beta)^2}{8 L_4^4}]$ when $\beta <L_2$.
 This shows that 
 \begin{align} \label{d2}
 D_2^\tau(\xi^q) \ge 0, \quad \xi \ge 0.  
 \end{align}
Combing the estimates \eqref{est-con}, \eqref{d1}, and \eqref{d2}, we conclude  
 \begin{align*} 
\Big(\frac{-L_2 \xi^q}{(1+\tau \xi^{2q})^{1/2}} 
+\frac{2L_4^2 \tau \xi^{2 q}}{1+\tau \xi^{2q}} \Big)\xi^2
\le  -\frac{\beta \alpha^q}{\sqrt{1+\alpha^{2q}}}\xi^2,
\quad \xi >\alpha> 0,  
 \end{align*}
 for all $\beta \in (0, L_2)$, $\tau \in (0, \frac{(L_2-\beta)^2}{8 L_4^4}]$,  and $\alpha>0$.

 For another part $\xi \in [0, \alpha]$, $D_1^\tau(\xi^q) \ge 0$ and  $D_2^\tau(\xi^q) \ge 0$, then we have
 \begin{align*}
 	\Big(\frac{-L_2 \xi^q}{(1+\tau \xi^{2q})^{1/2}} 
+\frac{2L_4^2 \tau \xi^{2 q}}{1+\tau \xi^{2q}} \Big)\xi^2
&=( -\beta D_1^\tau(\xi^q)  - \xi^q D_2^\tau(\xi^q))\xi^2\\
&\le 0= \Big(
-\frac{\beta \alpha^q}{\sqrt{1+\alpha^{2q}}}+
\frac{\beta \alpha^q}{\sqrt{1+\alpha^{2q}}}
\Big)\xi^2\\
&\le 
-\frac{\beta \alpha^q}{\sqrt{1+\alpha^{2q}}}\xi^2+\frac{\beta \alpha^{q+2}}{\sqrt{1+\alpha^{2q}}},
\quad 0 \le \xi \le \alpha.
 \end{align*}
 Therefore,
 \begin{align} \label{in-con}
 	\Big(\frac{-L_2 \xi^q}{(1+\tau \xi^{2q})^{1/2}} 
+\frac{2L_4^2 \tau \xi^{2 q}}{1+\tau \xi^{2q}} \Big)\xi^2
\le -\frac{\beta \alpha^q}{\sqrt{1+\alpha^{2q}}}\xi^2+\frac{\beta \alpha^{q+2}}{\sqrt{1+\alpha^{2q}}},
\quad \xi \ge 0,
 \end{align}
 from which we derive \eqref{coe-tau+} with $T_1:=L_1+ 2 L_3^2 \tau+\frac{\beta \alpha^{q+2}}{\sqrt{1+\alpha^{2q}}}$ and $T_2:=\frac{\beta \alpha^q}{\sqrt{1+\alpha^{2q}}}$. 
 \end{proof} 
 
 The following one-step Lyapunov estimate indicates that the original Lyapunov function $V: \rr^d \rightarrow [0,\infty)$ of Eq. \eqref{sde} defined by $V(x)=|x|^2$, $x \in \rr^d$, can be inherited by the TEM scheme \eqref{tem}, even though $V$ is not a Lyapunov function of the classical EM scheme as indicated in \cite{HJK11,MSH02}.
 We also note that the Lyapunov functions of implicit schemes, e.g., the stochastic theta scheme in \cite{LL25}, 
 take a form distinct from that of a square function.

 \begin{theorem} \label{tm-lya}
Let $X_0$ be $\FFF_0$-measurable such that $\mathbb{ E }|X_0|^2 < \infty$ and  the condition \eqref{coe-tau+} hold.
 Then there exists $\tau_{\max} \in (0, 1]$ such that for any $\tau \in (0, \tau_{\max}]$, there holds that 
 	\begin{align} \label{lya}
\mathbb{ E } [|Z_n|^2 | \FFF_{n-1}]
 	\le  (1- T_2 \tau) |Z_{n-1}|^2+ T_1 \tau, \quad n \in \nn_+.
 	\end{align}
 \end{theorem}
 
 \begin{proof}
For $n \in \nn_+$, taking square on both sides of the TEM scheme \eqref{tem}, we have 
 \begin{align*} 
|Z_n|^2&=|Z_{n-1}|^2+ 2\tau \<Z_{n-1},b_\tau(Z_{n-1})\>  
+ \tau^2 | b_\tau(Z_{n-1})|^2
+ |\sigma_\tau(Z_{n-1}) \delta_{n-1} W|^2  \\
& \quad + 2 \<Z_{n-1}+ b_\tau(Z_{n-1}) \tau, \sigma_\tau( Z_{n-1})\delta_{n-1} W\>.  
 \end{align*}
Assume that $\mathbb{ E }|Z_{n-1}|^2 < \infty$.
As $b_\tau$ and $\sigma_\tau$ are both have at most linear growth, we have $b_\tau(Z_{n-1}) \in L^2(\Omega; \rr^d)$ and $\sigma_\tau( Z_{n-1}) \in L^2(\Omega; \rr^{d \times m})$.
Note that $Z_{n-1}$ is $\FFF_{n-1}$-measurable and $\delta_{n-1} W$ is independent of $\FFF_{n-1}$, so 
 \begin{align*} 
 \ee[2 \<Z_{n-1}+ b_\tau(Z_{n-1}) \tau, \sigma_\tau( Z_{n-1})\delta_{n-1} W\>|\FFF_{n-1}]=0.  
 \end{align*}
Taking the conditional expectation  $\ee[\cdot |\FFF_{n-1}]$  on both sides of the above equation and using It\^o isometry, the condition \eqref{coe-tau+}, and the above equality, we obtain
 \begin{align*} 
 \ee[|Z_n|^2|\FFF_{n-1}]
& \le |Z_{n-1}|^2 + \tau [2\<Z_{n-1}, b_\tau(Z_{n-1})\>+\|\sigma_\tau(Z_{n-1})\|^2 
 + \tau | b_\tau(Z_{n-1})|^2] \\ 
&\le (1-T_2 \tau)|Z_{n-1}|^2 + T_1 \tau.
 \end{align*} 
 
This completes the proof of \eqref{lya} by  the condition $\mathbb{ E }|X_0|^2 < \infty$.
\end{proof}
 
 \begin{remark} \label{rk-lya} 
The estimate \eqref{lya} indicates that the natural Lyapunov function $V: \rr^d \to [0,\infty)$ of Eq. \eqref{sde} defined by
$V(x)=\|x\|^2$, $x \in \rr^d$, is also a Lyapunov function of the TEM scheme \eqref{tem}.

Moreover, if $0<T_2\tau \le 1$, then the following uniform estimate holds:
       \begin{align*}
\ee |Z_n|^2  
 	\le  T_1/T_2+e^{- T_2 \tau n} \ee |X_0|^2, \quad \forall~ n \in \nn_+.
      \end{align*}   
 \end{remark}

If the diffusion function $\sigma$ is of linear growth, then we do not need to tame $\sigma$ and consider the following drift-TEM scheme:
\begin{align} \label{tem-}
 	{\widetilde Z}_n={\widetilde Z}_{n-1}+ b_\tau ({\widetilde Z}_{n-1}) \tau
 	+\sigma({\widetilde Z}_{n-1})\delta_{n-1} W, 
	~ n \in \nn_+; \quad Z_0=X_0.
 \end{align} 
The above TEM scheme \eqref{tem-} has also been studied in \cite{NNZ25} for Eq. \eqref{sde} driven by additive noise.
Especially, the tamed drift function takes the form in \eqref{bs-tau} or \eqref{bs-tau+},
thus there exist positive constants $\widetilde T_1$ and $\widetilde T_2$ such that for all $x \in \rr^d$,
\begin{align}  \label{coe-tau+-linear}
2 \<x,b_\tau(x)\> + \|\sigma(x)\|^2 + \tau |b_\tau(x)|^2 
& \le \widetilde T_1 - \widetilde T_2 |x|^2. 
\end{align}

 Similar arguments in the proof of Theorem \ref{tm-lya} yield the following Lyapunov structure of the TEM scheme \eqref{tem-}.  
 
  \begin{corollary} \label{cor-lya}
Let $X_0$ be $\FFF_0$-measurable such that $\mathbb{ E }|X_0|^2 < \infty$ and the condition \eqref{coe-tau+-linear} hold.
 Then there exists $\tau_{\max} \in (0, 1]$ such that for any $\tau \in (0, \tau_{\max}]$, there holds that 
 	\begin{align} \label{lya-}
\mathbb{ E } [|\widetilde Z_n|^2 | \FFF_{n-1}]
 	\le  (1- \widetilde T_2 \tau) |\widetilde Z_{n-1}|^2
  + \widetilde T_1 \tau, \quad n \in \nn_+.
 	\end{align}
 \end{corollary}

\subsection{Geometric Ergodicity of TEM}

In this part, we focus on the geometric ergodicity of the TEM schemes \eqref{tem} and \eqref{tem-} under the following continuity and non-degeneracy conditions. 
  
\begin{assumption} \label{ap-non}
The functions $b_\tau$ and $\sigma_\tau$ are continuously differentiable, and for any $x \in \rr^d$, $\sigma_\tau(x) \sigma_\tau(x)^\top $ is positive definite on $\rr^d$.
\end{assumption}

As noted at the beginning of this section, $\{Z_n, \FFF_{t_n}\}_{n \in \nn}$ and $\{{\widetilde Z}_n, \FFF_{t_n}\}_{n \in \nn}$ generated by the TEM schemes \eqref{tem} and \eqref{tem-}, respectively, are two (time-homogeneous) Markov chains.  
Denote by $P: \rr^d \times \BB(\rr^d) \to [0, 1]$ and ${\widetilde P}: \rr^d \times \BB(\rr^d) \to [0, 1]$ the transition kernels of the Markov chain $\{Z_n, \FFF_{t_n}\}_{n \in \nn}$ and $\{{\widetilde Z}_n, \FFF_{t_n}\}_{n \in \nn}$, respetively.
 Then for any $n \in \nn$, 
 \begin{align} 
 P(x,A)=\mathbb{P}(Z_{n+1} \in A|Z_n=x),
 \quad \forall~ x \in \rr^d,
~ A \in \BB(\rr^d), \label{kernel} \\
{\widetilde P}(x,A)=\mathbb{P}({\widetilde Z}_{n+1} \in A|{\widetilde Z}_n=x),
 \quad \forall~ x \in \rr^d, \label{kernel-}
~ A \in \BB(\rr^d).
 \end{align}
 
 Let us first consider the TEM scheme \eqref{tem}. 
 For $ x \in \rr^d$ and $A \in \BB(\rr^d)$, 
\begin{align*}
P(x,A)&=\mathbb{P} (Z_{n+1} \in A |Z_n=x)
=\mu_{x+b_\tau(x),\sigma_\tau(x)\sigma_\tau(x)^\top \tau}(A),
 \end{align*}
 as $x+b_\tau(x)+\sigma_\tau(x) \delta_n W \sim \mathcal N(x+b_\tau(x),\sigma_\tau(x)\sigma_\tau(x)^\top \tau)$, where $\mathcal{N}(y, Q)$ and $\mu_{a,b}$ denote the normal distribution and Gaussian measure in $\rr^d$ with mean $y \in \rr^d$ and variance operator $Q \in \mathcal{L}(\rr^d)$, respectively. 
It is well-known that all non-degenerate $d$-dimensional Gaussian measures are equivalent.
Thus, $P$ is regular in the sense that all transition probabilities of \eqref{tem} are equivalent, and by Doob theorem, it possesses at most one invariant measure (see, e.g., \cite[Proposition 7.4]{Da06}). 
Similar uniqueness of the invariant measure, if it exists, holds for the TEM scheme \eqref{tem-}. 
 
In combination with the Lyapunov structure in Theorem \ref{tm-lya} and Corollary \ref{cor-lya}, and the proof of \cite[Corollary 3.2]{LL25}, we have the following geometric ergodicity of the TEM schemes \eqref{tem} and \eqref{tem-}, respectively.

\begin{theorem} \label{tm-erg}
Let Assumption \ref{ap-non} hold. 
\begin{enumerate}
\item
Under Assumption \ref{ap}, there exists ${\widetilde \tau}_{\max} \in (0, 1]$ such that for any $\tau \in (0, {\widetilde \tau}_{\max}]$, the TEM scheme \eqref{tem} is geometrically ergodic with respect to a unique invariant measure $\pi_\tau$ in $\BB(\rr^d)$, i.e., there exists $r_0 \in (0,1)$ and $\kappa \in (0,\infty)$, depending on $\tau$, such that for any measurable function $\phi: \rr^d \to \rr$ with $|\phi|\leq 1+|\cdot|^2$,
          \begin{align*}
          \abs{\ee \phi(Z_n^x)-\pi_\tau(\phi)}
          \leq \kappa r_0^n (1+|x|^2), \quad \forall~ x \in \rr^d.
          \end{align*}    

\item
Under the conditions \eqref{b-grow}, \eqref{b-coe}, and \eqref{s-grow}, there exists ${\widetilde \tau}_{\max}' \in (0, 1]$ such that for any $\tau \in (0, {\widetilde \tau}_{\max}']$, the TEM scheme \eqref{tem-} is geometrically ergodic.  
\qed
\end{enumerate} 
\end{theorem}

\begin{remark}
Assumption \ref{ap-non} is valid for the examples \eqref{bs-tau} and \eqref{bs-tau+} provided that $b$ and $\sigma$ are continuously differentiable and for any $x \in \rr^d$, $\sigma(x) \sigma(x)^\top $ is positive definite on $\rr^m$. 
\end{remark}

 \section{Error Estimates}
 \label{sec3}
 
In this section, we aim to establish quantitative error estimates between the TEM schemes  \eqref{tem} and \eqref{tem-} and Eq. \eqref{sde} in the $L^{2p}(\Omega)$-norm for $p \ge 2$ over the finite time interval $[0, T]$, where $T$ is a given positive number.
To this end, 
 We assume that the coefficients are monotone and coercive, and satisfy a polynomial growth condition in terms of the continuity modulus.

 \begin{assumption} \label{ap+}   
There exist positive constants $L_0'$, $L_7$, $L_8$, $L_9$, $L_{10}$, $q$, and $p^*>p \ge 2$ such that for all $x,y \in \rr^d$, 
\begin{align} 
& 2 \<x-y,b(x)-b(y)\>+ (2 p^*-1) \|\sigma(x)-\sigma(y)\|^2 
\le L_0' |x-y|^2, \label{mon} \\
& 2 \<x,b(x)\> + (2p-1) \|\sigma(x)\|^2 
\le L_7 + L_8 |x|^2, \label{coe-} \\ 
& |b(x)-b(y)| \le L_9(1+|x|+|y|)^q |x-y|,\label{b-lip}  \\
& \|\sigma(x)-\sigma(y)\|^2 \le L_{10} (1+|x|+|y|)^q |x-y|^2. \label{s-lip} 
	\end{align} 
\end{assumption}
 
 \begin{remark}\label{rk-ap+}
 \begin{enumerate}   
\item
The coupled monotone condition \eqref{mon} includes the drift and diffusion functions given in \eqref{bs-ex}  when $d=1$.
For example, let $2=p<p^*<3.5$, $b(x)=(1-x^2) x$, and $\sigma(x)=(1+x^2)/2$, $x \in \rr$.
Then
\begin{align*} 
& 2 \<x-y,b(x)-b(y)\>+ (2 p^*-1) \|\sigma(x)-\sigma(y)\|^2 \\
& = |x-y|^2 [2-2(x^2+xy+y^2)+ \frac{2 p^*-1}4 (x+y)^2] \\
& \le 2 |x-y|^2, \quad \forall~ x,y \in \rr, 
	\end{align*}
	which shows \eqref{mon} with $L_0'=2$.

\item 
Under the condition \eqref{coe}, Young inequality implies the coercive condition \eqref{coe-} with some $L_7$ and $L_8$ depending on $L_1$, $L_2$, $p$, and $q$.

 \item
The polynomial growth conditions \eqref{b-grow} and \eqref{s-grow+} follow immediately from \eqref{b-lip} and \eqref{s-lip}, respectively; they are equivalent, respectively, under the continuous differentiability condition of $b$ and $\sigma$.

\item 
Under the conditions \eqref{mon} and \eqref{coe-}, Eq. \eqref{sde} possesses
a unique $\mathbb{F}$-adapted solution $\{X(t)\}_{t\in[0,T]}$ with
continuous sample paths (see \cite{LR15}).
Moreover, if $\ee |X_0|^{2p} < \infty$ for some $p\ge 1$, then there exists a constant $C$ such that (see \cite{HMS02})
\begin{align}\label{bound-x}
 \ee|X(t)|^{2p}\le e^{Ct}   (Ct+\ee|X_0|^{2p}), \quad t \ge 0.
\end{align}
Here and after, $C$ denotes a generic constant independent of various discrete parameters and the terminal time $T$ that may differ in each appearance. 
Furthermore, if \eqref{b-lip} and \eqref{s-lip} hold, one has the following H{\"o}lder continuity of the exact solution:
(see \cite{Liu22}):
\begin{align}\label{hol-x}
 \ee|X(t)-X(s)|^{2p}\le e^{CT}   (1+CT+\ee|X_0|^{p(q+2)})(t-s)^p,
 \quad t, s \in [0, T].
\end{align}
\end{enumerate}
\end{remark}

To utilize the tools of stochastic analysis, we first define the continuous-time interpolation process $\{X_\tau(t): t \ge 0\}$ of the TEM scheme \eqref{tem} by
\begin{align} \label{tem-con}
	{\rm d }X_\tau ( t )
	=b_\tau (X_\tau ( \kappa_n(t)) ){\rm d } t 
	+ \sigma_\tau (X_\tau ( \kappa_n(t)) ){\rm d }W ( t ), ~ t \ge 0;
	\quad X_\tau(0)=X_0,
\end{align}
where $\kappa_n(t)= t_n:=\tau n$ if $t\in [t_n ,t_{n+1})$.
Obviously, $X_\tau$ is an It\^o process and $Z_n=X_\tau(t_n)$ for all $n \in \nn$.

To derive moment estimates, H\"older continuity, and strong error bounds, we introduce the following assumption. 

\begin{assumption}\label{ap*}

There exist positive constants ${\widetilde L}_i$, $i=1,2,\cdots,8$, and $r$, independent of $\tau$,
 such that for any $x \in \rr^d$ and $\tau \in (0, 1)$,  
\begin{align}
2 \<x, b_\tau(x)\> + (2p-1) \|\sigma_\tau(x)\|^2 
& \le {\widetilde L}_1 + {\widetilde L}_2 |x|^2, \label{coe-tau} \\
\tau|b_\tau(x)|^2  
& \le {\widetilde L}_3+{\widetilde L}_4 |x|^2,   \label{btau-grow} \\
\sqrt{\tau}\|\sigma_\tau(x)\|^2
& \le \widetilde L_5+\widetilde L_6|x|^2,  \label{stau-grow} \\
|b(x)-b_\tau(x)|^2
& \le {\widetilde L}_7 \tau^2(1+|x|^r)|b(x)|^2 , \label{b-btau} \\
\|\sigma(x)-\sigma_\tau(x)\|^2
& \le {\widetilde L}_8 \tau(1+|x|^r)\|\sigma(x)\|^2.  \label{s-stau}
\end{align} 
\end{assumption}

\begin{remark}\label{rk-ap*}
Under Assumption \ref{ap+}, the above Assumption \ref{ap*} is valid for the tamed functions $b_\tau$ and $\sigma_\tau$ given in \eqref{bs-tau} or \eqref{bs-tau+}.
The first three conditions in Assumption \ref{ap+} are obvious.
To show the other two conditions, we define
\begin{align*}
g(\xi):=(1-(1+\xi)^{-1/2}) ^2, \quad \xi \ge 0. 
\end{align*}
Direct computations yield   
\begin{align*}
g(\xi) & =\frac{((1+\xi)^{1/2}-1)^2((1+\xi)^{1/2}+1)^2}{(1+\xi)((1+\xi)^{1/2}+1)^2} 
 =\frac{\xi^2}{(1+\xi)((1+\xi)^{1/2}+1)^2} 
 \le \frac14 \xi^2,  
\end{align*}
for any $\xi \ge 0$.  
Then we use the representation \eqref{bs-tau} to derive 
\begin{align*}
|b(x)-b_\tau(x)|^2
=g(\tau|x|^{2q}) |b(x)|^2
\le \frac14 \tau^2|x|^{4q}|b(x)|^2 ,
\end{align*}
and
\begin{align*}
\|\sigma(x)-\sigma_\tau(x)\|^2
=\sigma(\sqrt{\tau}|x|^q) \|\sigma(x)\|^2
\le \frac14 \tau|x|^{2q}\|\sigma(x)\|^2 ,
\end{align*} 
for any $x \in \rr^d$, which show \eqref{b-btau} and \eqref{s-stau} (with ${\widetilde L}_7={\widetilde L}_8=1/4$ and $r=4q$).

The treatment of Example \eqref{bs-tau+} is identical,  and we omit the details. 
\end{remark}

\subsection{Moment's Estimate and H\"older Continuity}

We begin with the following moment's estimate of the interpolation process $X_\tau$ defined in \eqref{tem-con}.

In combination with the conditions \eqref{coe-} and \eqref{mon}, it is clear that Eq. \eqref{sde} is well-posed and possesses bounded algebraic moments provided the initial datum has bounded algebraic moments, see, e.g., \cite[Lemma 3.2]{HMS02} and \cite[Proposition 3.1]{Liu22}.

\begin{lemma}
	\label{lm-bound-xtau}
Let $X_0$ be $\FFF_0$-measurable such that $X_0 \in L^{2p}(\Omega; \rr^d)$ with $p \ge 2$, and \eqref{b-grow} and \eqref{coe-} hold.
Then there exists a positive constant $C$ such that  
	\begin{align} \label{bound-xtau}
\ee |X_\tau(t)|^{2p} \le e^{Ct}  (C t+ \ee |X_0|^{2p}), \quad t \ge 0.
	\end{align} 
\end{lemma}

\begin{proof}
Let $t \ge 0$.
Applying It\^o formula to the functional $|X_\tau|^{2p}$, we get 
	\begin{align*}
|X_\tau(t)|^{2 p}&= |X_0|^{2 p}+2p \int_0^t  |X_\tau|^{2p-2}\< X_\tau, b_\tau(X_\tau(\kappa_n))\> {\rm d } s \\
&\quad+2p \int_0^t  |X_\tau|^{2p-2} \<X_\tau, \sigma_\tau(X_\tau(\kappa_n)) {\rm d } W\> \\
&\quad+p\int_0^t |X_\tau|^{2p-2}\|\sigma_\tau(X_\tau(\kappa_n))\|^2{\rm d}s\\
&\quad+ 2p(p-1) \int_0^t  |X_\tau|^{2(p-2)}|X_\tau^\top
\sigma_\tau(X_\tau(\kappa_n))|^2{\rm d } s, 
	\end{align*}
where we omit the integration variable here and after when there is an integration to lighten the notation.
Taking expectations on both sides of the above equation and using \eqref{coe-tau}, the martingale property of the It\^o integral above, and Young inequality, we obtain
\begin{align*}
\ee |X_\tau(t)|^{2p}
&\le \ee |X_0|^{2p }+ 2p\ee \int_0^t  |X_\tau|^{2p-2}\<X_\tau-X_\tau(\kappa_n), b_\tau(X_\tau(\kappa_n))\> {\rm d } s \\
&\quad + p\ee\int_0^t  |X_\tau|^{2p-2} [2\<X_\tau(\kappa_n), b_\tau(X_\tau(\kappa_n))\> +(2p-1) \|\sigma_\tau(X_\tau(\kappa_n))\|^2] {\rm d } s\\ 
& \le \ee |X_0|^{2p }+ 2p\ee \int_0^t  |X_\tau|^{2p-2}\<X_\tau-X_\tau(\kappa_n), b_\tau(X_\tau(\kappa_n))\> {\rm d } s \\
&\quad + p\ee\int_0^t  |X_\tau|^{2p-2} ({\widetilde L}_1 +{\widetilde L}_2 |X_\tau(\kappa_n)|^2) {\rm d } s\\ 
&\le Ct+ \ee |X_0|^{2p}+C \int_0^t \sup_{ u \le  s}\ee|X_\tau (u)|^{2p}{\rm d } s + 2p J,
\end{align*}
where $J:=\ee \int_0^t  |X_\tau|^{2 (p-1)}\< X_\tau-X_\tau(\kappa_n), b_\tau(X_\tau(\kappa_n))\> {\rm d } s$.

To estimate the term $J$, we use  \eqref{tem-con} to rewrite
\begin{align*}
J&=\ee\int_0^t |X_\tau|^{2 (p-1)}
|b_\tau(X_\tau(\kappa_n))|^2 (s-\kappa_n) {\rm d } s\\
&\quad + \ee\int_0^t |X_\tau|^{2 (p-1)} b_\tau(X_\tau(\kappa_n)) \sigma_\tau(X_\tau(\kappa_n))(W-W(\kappa_n)) {\rm d} s
=:J_1+J_2.
\end{align*}
Then we use \eqref{btau-grow} and Young inequality to bound $J_1$:
\begin{align*}
J_1 & \le \ee\int_0^t |X_\tau|^{2p-2}(\tau|b_\tau(X_\tau(\kappa_n))|^2){\rm d} s\\ 
&\le \ee\int_0^t |X_\tau|^{2 (p-1)}({\widetilde L}_4 |X_\tau(\kappa_n)|^2 +{\widetilde L}_3) {\rm d } s\\ 
&\le Ct+C \int_0^t \sup_{ u \le  s}\ee|X_\tau (u)|^{2p}{\rm d } s.
\end{align*}
For the term $J_2$, the fact the centered Brownian increment $W(s)-W(\kappa_n(s))$ is independent of $X_\tau(\kappa_n(s))$ for all $s \in (0, t)$ and It\^o formula applied to the functional $|X_\tau|^{2p-2}$ yield that  
\begin{align*}
J_2&=2p\ee\int_0^t  (|X_\tau|^{2p-2}-|X_\tau(\kappa_n)|^{2 (p-1)}) b_\tau(X_\tau(\kappa_n)) \sigma_\tau(X_\tau(\kappa_n))(W-W(\kappa_n)){\rm d } s\\ 
&=2p \ee\int_0^t (2(p-1)\int_{\kappa_n(s)}^s|X_\tau|^{2(p-2)}\<
X_\tau,b_\tau(X_\tau(\kappa_n))\>{\rm d } v\\
&\quad+2(p-1) \int_{\kappa_n(s)}^s|X_\tau|^{2(p-2)}\<X_\tau,\sigma_\tau(X_\tau(\kappa_n)){\rm d } W\>\\
&\quad+ (p-1) \int_{\kappa_n(s)}^s|X_\tau|^{2(p-2)}\|\sigma_\tau (X_\tau(\kappa_n))\|^2{\rm d } v\\
&\quad+ 2(p-1)(p-2) \int_{\kappa_n(s)}^s|X_\tau|^{2(p-3)}|X_\tau^\top\sigma_\tau(X_\tau(\kappa_n))|^2{\rm d } v)\\
&\quad \times b_\tau(X_\tau(\kappa_n)) \sigma_\tau(X_\tau(\kappa_n))(W-W(\kappa_n)) {\rm d } s\\
&=:J_{21} +J_{22}+J_{23}+J_{24}.
\end{align*}
Following we estimate $J_{21}$ to $J_{24}$ by using the estimates \eqref{btau-grow}, \eqref{stau-grow},
and H\"older and Young inequalities.
More precisely,
\begin{align*}
J_{21}
&\le C \ee\int_0^t \Big(\int_{\kappa_n(s)}^s|X_\tau|^{2p-3}
{\rm d } v\Big) |b_\tau(X_\tau(\kappa_n))|^2
|\sigma_\tau(X_\tau(\kappa_n))(W-W(\kappa_n))| {\rm d } s\\
&\le C \int_0^t  \tau^{-1} [\ee (\tau |b_\tau(X_\tau(\kappa_n))|^2 \int_{\kappa_n(s)}^s|X_\tau|^{2p-3} {\rm d } v )^{\frac{2p}{2p-1}}]^{\frac{2p-1}{2p}} \\
&\qquad\qquad   \times [\ee (\tau \|\sigma_\tau(X_\tau(\kappa_n))\|^2)^p]^{\frac{1}{2p}}  {\rm d } s\\
&\le C \ee\int_0^t [\tau |b_\tau(X_\tau(\kappa_n))|^2]^{\frac{2p}{2p-1}} 
\Big(\tau^{-1}\int_{\kappa_n(s)}^s|X_\tau|^{2p-3}{\rm d }v \Big)^{\frac{2p}{2p-1}}{\rm d } s \\
&\quad + C \int_0^t  \ee [\tau \|\sigma_\tau(X_\tau(\kappa_n))\|^2]^p {\rm d } s\\
&\le C \ee\int_0^t (\tau |b_\tau(X_\tau(\kappa_n))|^2)^p 
+ (\sqrt{\tau} \|\sigma_\tau(X_\tau(\kappa_n))\|^2)^p {\rm d } s \\ 
& \quad + \ee\int_0^t (\tau^{-1}\int_{\kappa_n(s)}^s|X_\tau|^{2p-3}{\rm d }v )^{\frac{2p}{2p-3}}{\rm d } s \\ 
&\le Ct+C\int_0^t \sup_{ u \le  s}\ee|X_\tau|^{2p}{\rm d } s.
\end{align*}
Similarly, the It\^o isometry implies that 
\begin{align*}
J_{22}
&\le C \ee\int_0^t (\tau^{-1}\int_{\kappa_n(s)}^s|X_\tau|^{2p-3}{\rm d } v) \sqrt{\tau}\| \sigma_\tau(X_\tau(\kappa_n))\|^2|\sqrt \tau b_\tau(X_\tau(\kappa_n))|{\rm d } s\\
&\le C \ee \int_0^t (\tau^{-1} \int_{\kappa_n(s)}^s
|X_\tau|^{2p-3}{\rm d } v)^{\frac{2p}{2p-3}} \\
&\quad + C \ee \int_0^t(\sqrt{\tau}\|\sigma_\tau(X_\tau(\kappa_n))\|^2
|\sqrt{\tau} b_\tau(X_\tau(\kappa_n))| )^{\frac{2p}{3}} {\rm d } s\\
&\le Ct+C\int_0^t \sup_{ u \le  s}\ee|X_\tau|^{2p}{\rm d } s	.
\end{align*}
Finally, when $p>2$, 
\begin{align*}
J_{23}+J_{24}
&\le C\int_0^t (\ee(\sqrt{\tau}\|\sigma_\tau(X_\tau(\kappa_n))\|^2
\sqrt \tau|b_\tau(X_\tau(\kappa_n))| \\
&\quad\times \tau^{-1} \int_{\kappa_n(s)}^s
|X_\tau|^{2(p-2)}{\rm d } v )^{\frac{2p}{2p-1}})^{\frac{2p-1}{2p}}
(\tau^p\ee\|\sigma_\tau(X_\tau(\kappa_n))\|^{2p})^{\frac{1}{2p}}{\rm d } s
\\ 
&\le C\ee\int_0^t (\sqrt{\tau}\|\sigma_\tau(X_\tau(\kappa_n))\|^2
\sqrt \tau|b_\tau(X_\tau(\kappa_n))|)^{\frac{2p}{2p-1}} \\
&\quad\times(\tau^{-1}\int_{\kappa_n(s)}^s|X_\tau|^{2(p-2)}
{\rm d } v)^{\frac{2p}{2p-1}}{\rm d } s 
+ C\tau^p\int_0^t \ee \|\sigma_\tau(X_\tau(\kappa_n))\|^{2p} {\rm d } s\\
&\le C\ee\int_0^t (\sqrt{\tau}\|\sigma_\tau(X_\tau(\kappa_n))\|^2
\sqrt{\tau}|b_\tau(X_\tau(\kappa_n))|)^{\frac{2p}{3}}
+ (\tau \|\sigma_\tau(X_\tau(\kappa_n))\|^2)^p {\rm d } s \\
& \quad + C\ee\int_0^t (\tau^{-1}\int_{\kappa_n(s)}^s|X_\tau|^{2(p-2)}
{\rm d } v)^{\frac{p}{p-2}}{\rm d } s  \\
&\le Ct+C\int_0^t \sup_{ u \le  s}\ee|X_\tau|^{2p}{\rm d } s	.
\end{align*}

Note that $J_{24}$ vanishes when $p=2$.
Combining the above estimates, we obtain
\begin{align*}
\ee|X_\tau(t)|^{2p}&\le Ct+\ee|X_0|^{2p}+C\int_0^t \sup_{u\le  s}
\ee|X_\tau|^{2p}{\rm d } s	.
\end{align*}
Finally, the application of Gronwall lemma yields \eqref{bound-xtau}.
\end{proof}

Next, in any subinterval, we have the following H\"older continuity of the interpolation process \eqref{tem-con}.

\begin{lemma} \label{lm-hol-xtau}
Let $X_0$ be $\FFF_0$-measurable such that $X_0 \in L^{p(q+2)}(\Omega; \rr^d)$ with $p \ge 2$, and \eqref{b-grow} and \eqref{coe-} hold.
Then there exists a positive constant $C$ such that 
	\begin{align} \label{hol-xtau}
\ee |X_\tau(t)-X_\tau ( \kappa_n(t))|^{2p} 
\le e^{Ct} (1+Ct+ \ee |X_0|^{p(q+2)}) \tau^{p}, \quad n \in \nn_+.
	\end{align} 
	Consequently, for any $s, t \in [t_n, t_{n+1})$ with some $n \in \nn$, 
	\begin{align} \label{hol-xtau+}
\ee |X_\tau(t)-X_\tau (s)|^{2p} 
\le e^{C(t \vee s)} (1+C (t \vee s)+ \ee |X_0|^{p(q+2)}) \tau^{p}.
	\end{align} 
\end{lemma}

\begin{proof}
By \eqref{tem-con} and the elementary inequality $(a+b)^{2p} \le 2^{2p-1}(a^{2p}+b^{2p})$ for $a, b \ge 0$, we have 
\begin{align*}
&\ee \Big|X_\tau(t)-X_\tau (\kappa_n(t))\Big|^{2p }\\
&=\ee \Big|\int_{\kappa_n(t)}^{t}b_\tau( X_\tau (\kappa_n) ){\rm d } s + 
\int_{\kappa_n(t)}^{t}\sigma( X_\tau (\kappa_n) ){\rm d }W ( s )\Big|^{2p}\\
&\le C\ee\Big|\int_{\kappa_n(t)}^{t}b_\tau( X_\tau (\kappa_n) ){\rm d } s \Big|^{2p}+ C\ee\Big|\int_{\kappa_n(t)}^{t}\sigma( X_\tau (\kappa_n) ){\rm d }W ( s )\Big|^{2p}.
\end{align*} 
According to H\"older and BDG inequalities and the conditions \eqref{btau-grow} and \eqref{s-grow+} followed by the conditions \eqref{b-grow} and \eqref{coe-}, we have
	\begin{align*}
&\ee |X_\tau(t)-X_\tau (\kappa_n(t))|^{2p }\\
&\le C\tau^{2p-1}\ee\int_{\kappa_n(t)}^{t}|b_\tau(X_\tau (\kappa_n))|^{2p}
{\rm d }s + C\tau^{p-1}\ee\int_{\kappa_n(t)}^{t}\|\sigma_\tau (X_\tau (\kappa_n))\|
^{2p}{\rm d } s \\
&\le C\tau^{p-1}\ee\int_{\kappa_n(t)}^{t}({\widetilde L}_4|X_\tau( \kappa_n)|^2+{\widetilde L}_3)^p + (L_6|X_\tau( \kappa_n) |^{q+2}+L_5)^p{\rm d } s \\
&\le C\tau^p + C \tau^p \ee |X_\tau( \kappa_n)|^{p(q+2)},
\end{align*} 
which shows \eqref{hol-xtau} using the moment's estimate \eqref{bound-xtau}.  

Fanally, for any $s, t \in [t_n, t_{n+1})$ with some $n \in \nn$, it is clear that $\kappa_n(t)=\kappa_n(s)=t_n$. 
Consequently,  
	\begin{align*}
\ee |X_\tau(t)-X_\tau (s)|^{2p} 
& \le C |X_\tau(t)-X_\tau (\kappa_n(t))|^{2p}+C |X_\tau(s)-X_\tau (\kappa_n(s))|^{2p}.
	\end{align*} 
This shows \eqref{hol-xtau+} using the H\"older continuity estimate \eqref{hol-xtau}.  
\end{proof}

\subsection{Strong Error Estimates of Multiplicative Noise Case}

 In this part, we establish a finite-time strong error estimate for the convergence of the interpolation process \eqref{tem-con} towards Eq. \eqref{sde}.

 \begin{theorem}  \label{tm-err}
 Let  $X_0$ be $\FFF_0$-measurable with
 $X_0 \in L^{pr_1^*}(\Omega; \rr^d)$ for $p \ge 2$ and 
 $r_1^*=\max\{2q+2+r,3q+2\}$, and Assumptions \ref{ap+}, \ref{ap*} hold.
Then there exists a positive constant $C$ such that 
 	\begin{align} \label{err}
\ee |X({t_n})-X_\tau({t_n})|^{2p}
 \le  C e^{C t_n}\tau^p (1+ C t_n+ \ee |X_0|^{pr_1^*}), \quad n \in \nn_+.
 	\end{align}   
 \end{theorem}

 \begin{proof}
 Using It\^o formula to the functional $|X-X_\tau|^{2 p}$,  
we have
\begin{align} \label{ito-err}
& |X({t_n})-X_\tau({t_n})|^{2p} \nonumber  \\
&=2p \int_{0}^{t_n}|X-X_\tau|^{2p-2}
\langle X-X_\tau,b(X)-b_\tau(X_\tau( \kappa_n)) \rangle {\rm d } s \nonumber  \\
&\quad+2p \int_{0}^{t_n}|X-X_\tau|^{2p-2}\langle X-X_\tau,
(\sigma(X)-\sigma_\tau (X_\tau( \kappa_n))) {\rm d } W\rangle \nonumber  \\
&\quad+p\int_{	0}^{t_n}|X-X_\tau|^{2p-2}\|\sigma (X)-\sigma_\tau (X_\tau( \kappa_n)) \|^2 {\rm d } s  \nonumber  \\
&\quad+2p(p-1)\int_{	0}^{t_n}|X-X_\tau|^{2(p-2)}|(X-X_\tau)^\top
(\sigma (X)-\sigma_\tau (X_\tau( \kappa_n)) )|^2 {\rm d } s.
\end{align}
Taking expectations on both sides of the above equation and using the elementary inequality $2ab \le \epsilon a^2+ \frac{1}{\epsilon}b^2$ for $a, b \in \rr$( where here and after $\epsilon$ denotes an arbitrary infinitesimal positive constant),  together with the martingale property of the It\^o integral above, we obtain
\begin{align*}
& \ee|X({t_n})-X_\tau({t_n})|^{2p} \\
&\le 2p\ee \int_{0}^{t_n}
|X-X_\tau|^{2p-2} \langle X-X_\tau,
b(X)-b_\tau(X_\tau( \kappa_n)) \rangle {\rm d } s\\
&\quad+p(2p-1)\ee \int_{0}^{t_n}|X-X_\tau|^{2p-2}
\| \sigma (X)-\sigma_\tau(X_\tau( \kappa_n)) \|^2 {\rm d } s\\
 &\le p\ee \int_{0}^{t_n}|X-X_\tau|^{2p-2}
 [2 \langle X-X_\tau,b(X)-b(X_\tau) \rangle \\
& \qquad   \qquad   \qquad   \qquad   \qquad   \quad  
+ (2p-1)(1+\epsilon) \| \sigma (X)-\sigma (X_\tau) \|^2] {\rm d } s\\
&\quad+2p\ee \int_{	0}^{t_n}|X-X_\tau|^{2p-2}
\langle X-X_\tau,b(X_\tau)-b(X_\tau( \kappa_n)) \rangle {\rm d } s\\
&\quad+2p\ee \int_{	0}^{t_n}|X-X_\tau|^{2p-2}
\langle X-X_\tau,b(X_\tau( \kappa_n))-b_\tau(X_\tau( \kappa_n)) \rangle {\rm d } s\\
&\quad+2p(2p-1)(1+\frac{1}{\epsilon})\ee \int_{0}^{t_n}
|X-X_\tau|^{2p-2}\|\sigma (X_\tau)-\sigma (X_\tau( \kappa_n))\|^2 {\rm d } s\\
&\quad+2p(2p-1)(1+\frac{1}{\epsilon})\ee \int_{0}^{t_n}|X-
X_\tau|^{2p-2}\| \sigma (X_\tau( \kappa_n))
-\sigma_\tau(X_\tau( \kappa_n)) \|^2 {\rm d } s.
 	\end{align*}
It is clear that there exists  $\epsilon$ such that $(2p-1)(1+\epsilon)\le (2p^*-1)$, since $p<p^*$ by assumption.
Then we derive from the Assumptions \ref{ap+} and \ref{ap*} that 
\begin{align*}
& \ee|X({t_n})-X_\tau({t_n})|^{2p} \\ 
&\le C \ee \int_{0}^{t_n}|X-X_\tau|^{2p}{\rm d } s
+C\ee \int_{0}^{t_n}|b(X_\tau)- b(X_\tau(\kappa_n))|^{2p} {\rm d } s\\
 &\quad+C\ee \int_{0}^{t_n}|b(X_\tau(\kappa_n))
 -b_\tau(X_\tau(\kappa_n))|^{2p} {\rm d } s\\
&\quad+ C\ee \int_{0}^{t_n}\|\sigma(X_\tau)-\sigma(X_\tau( \kappa_n))\|^{2p} {\rm d } s\\
&\quad+ C\ee \int_{0}^{t_n} \|\sigma(X_\tau(\kappa_n))
 -\sigma_\tau(X_\tau( \kappa_n))\|^{2p}{\rm d } s \\
 &\le C \ee \int_{0}^{t_n}|X-X_\tau|^{2p}{\rm d } s\\
 &\quad+C\ee \int_{0}^{t_n}|X_\tau-X_\tau(\kappa_n)|^{2p}(1+|X_\tau|+
 |X_\tau(\kappa_n)|)^{2pq}{\rm d } s\\
&\quad+C\tau^{2p} \ee \int_{0}^{t_n}(1+|X_\tau(\kappa_n)|^r)^p 
(L_4|X_\tau(\kappa_n)|^{q+1}+L_3)^{2p}  {\rm d } s\\
&\quad+C\ee \int_{0}^{t_n}|X_\tau-X_\tau(\kappa_n)|^{2p}(1+|X_\tau|+
|X_\tau(\kappa_n)|)^{pq} {\rm d } s\\
&\quad+C\tau^p\ee \int_{0}^{t_n}(1+|X_\tau(\kappa_n)|^r)^p 
 (L_6|X_\tau(\kappa_n)|^{q+2}+L_5)^{p}{\rm d } s.
 \end{align*}
 By  the estimates \eqref{hol-xtau} and \eqref{bound-xtau}, we obtain
 \begin{align*}
& \ee|X({t_n})-X_\tau({t_n})|^{2p} \\
&\le C \ee \int_{0}^{t_n}|X-X_\tau|^{2p} {\rm d } s
+ C \tau^p \int_{0}^{t_n} (1+ \ee |X_\tau(\kappa_n)|^{p(2q+2+r)})  {\rm d } s \\
&\quad+C \int_{	0}^{t_n}\Big(\ee|X_\tau-X_\tau(\kappa_n)|^{\frac{2p(3q+2)}{q+2}}\Big)^{\frac{q+2}{3q+2}}\big(\ee(1+|X_\tau|+|X_\tau(\kappa_n)|)
^{p(3q+2)}\big)^{\frac{2q}{3q+2}}{\rm d } s\\  
&\le C \ee \int_{0}^{t_n}|X-X_\tau|^{2p} {\rm d } s
 +C e^{Ct_n}\tau^p(1+C t_n+ \ee |X_0|^{p(2q+2+r)}+\ee |X_0|^{p(3q+2)}).  
 	\end{align*}
An application of Gronwall inequality yields the strong error estimate \eqref{err}.
 \end{proof}

     \subsection{Strong Error Estimates of Additive Noise Case}

Now we consider the additive noise case, i.e., the diffusion coefficient of Eq. \eqref{sde} is a constant matrix ${\rm B} \in \rr^{d \times m}$: 
 \begin{align}\label{sde*}
    \dd X(t) = b(X(t))\dd t + {\rm B} \dd W(t), ~ t \ge 0.  
\end{align}
 In this case, we consider the following version of the drift-TEM scheme (i.e., \eqref{tem-} with $\sigma(x) ={\rm B}$, $x \in \rr^d$) applied to Eq. \eqref{sde*}:
 \begin{align} \label{tem*}
 	Z_n=Z_{n-1}+
 	\tau b_\tau (Z_{n-1})
 	+{\rm B} \delta_{n-1} W,
~ n \in \nn_+; \quad Z_0=X_0, 
 \end{align}
 where $b_\tau$ is defined as \eqref{bs-tau} or \eqref{bs-tau+}, as in the TEM scheme \eqref{tem}.
 Then the corresponding continuous-time interpolation of drift-TEM \eqref{tem-} can be written as 
 \begin{align} \label{tem-con*}
{\rm d}X_\tau(t)=b_\tau (X_\tau(\kappa_n(t))) {\rm d}t+ {\rm B} {\rm d}W(t), ~ t \ge 0; \quad X_\tau(0)=X_0.
 \end{align} 

In the additive noise case, we obtain a high-order convergence result under additional regularity and polynomial growth assumptions on the drift function $b$:
 \begin{align}\label{b-add-poly}
|\nabla^2 b(x) (y, z)| \le (L_4' |x|^{q-1}+L_3') |y| |z|,
\quad x, y, z \in \rr^d,
  \end{align}
  for some positive constants $L_3'$ and $L_4'$, where $\nabla^2$ denotes the Hessian matrix of $b$.  
  
For a twice differentiable function $b:\rr^d \rightarrow\rr^d$, Taylor formula combined with \eqref{tem-con*} implies 
\begin{align}\label{b-taylor}
& b(X_\tau(s))-b(X_\tau(\kappa_n(s))) \nonumber\\
& =\nabla b(X_\tau(\kappa_n(s))) {\rm B} (W(s)-W(\kappa_n(s))) \nonumber \\
& \quad + (s-\kappa_n) \nabla b(X_\tau(\kappa_n(s))) b_\tau(X_\tau(\kappa_n(s)))+R(X_\tau(s), X_\tau(\kappa_n(s))),
\end{align} 
where 
\begin{align}\label{R}
R(x,y)=\int_{0}^{1} (1-r)\nabla^2 b(y+r(x-y)) (x-y, x-y){\rm d}r,
\quad x, y \in \rr^d.
\end{align}
We first give an error estimate of the remainder $R$ defined in \eqref{R}.

\begin{lemma} 
Let  $X_0$ be $\FFF_0$-measurable with
$X_0 \in L^{2p(2q+1)}(\Omega; \rr^d)$ for $p \ge 2$ and $q>1$, and Assumptions \ref{ap+}, \ref{ap*} and the condition \eqref{b-add-poly} hold.
Then there exists a positive constant $C$ such that 
\begin{align} \label{err-R}
\ee|R(X_\tau,X_\tau(\kappa_n))|^{2p}
\le C e^{Ct_n}\tau^{2p}(1+C t_n+\ee|X_0|^{2p(2q+1)}), \quad n \in \nn_+.
\end{align}
\end{lemma}

\begin{proof} 
Using the representation \eqref{R}, the condition \eqref{b-add-poly}, H\"older inequality, and the estimates \eqref{bound-xtau} and \eqref{hol-xtau}, we obtain
\begin{align*}
& \ee|R(X_\tau(s), X_\tau(\kappa_n(s)))|^{2p} \\
&\le \ee\int_{0}^{1} |\nabla^2 b(X_\tau(\kappa_n)+
r(X_\tau-X_\tau(\kappa_n))(X_\tau-X_\tau(\kappa_n), 
X_\tau-X_\tau(\kappa_n))|^{2p} {\rm d}r\\ 
&\le \int_{0}^{1}(\ee[|X_\tau(\kappa_n)|^{q-1}+|X_\tau|^{q-1}+1]^{\frac{2p(2q+1)}{q-1}})^{\frac{q-1}{2q+1}} (\ee|X_\tau-X_\tau(\kappa_n)|^{\frac{4p(2q+1)}{q+2}})^{\frac{q+2}{2q+1}} {\rm d}r\\
&\le C e^{Ct_n}\tau^{2p}(1+C t_n+ \ee|X_0|^{2p(2q+1)}).
\end{align*}
This shows \eqref{err-R}.
\end{proof}

\begin{theorem} \label{tm-err+}
 Let  $X_0$ be $\FFF_0$-measurable with $X_0 \in L^{2p r_2^*}(\Omega; \rr^d)$ for $p \ge 2$ and $r_2^*=\max\{q+1+r/2, 2q+1\}$ with $q>1$ and $r$ given in \eqref{b-btau} and \eqref{s-stau}, and  the conditions in Theorem \ref{tm-err} and \eqref{b-add-poly} hold. There exists a positive constant $C$  such that
	\begin{align} \label{err-L2*}
		\ee |X({t_n})-X_\tau({t_n})|^{2p} 
		\le C e^{Ct_n } \tau^{2p} (1+C t_n+ \ee |X_0|^{2p r_2^*}), \quad n \in \nn_+.
	\end{align}  
\end{theorem}

\begin{proof}
Applying It\^o formula as in \eqref{ito-err}, we have
\begin{align}\label{ito-err*}
&\ee|X({t_n})-X_\tau({t_n})|^{2p}  \nonumber\\
&=2p\ee\int_{0}^{t_n}|X-X_\tau|^{2p-2}
\< X-X_\tau,b(X)-b_\tau(X_\tau( \kappa_n)) \> {\rm d } s \nonumber\\
&=2p\ee\int_{0}^{t_n}|X-X_\tau|^{2p-2}
\< X-X_\tau,b(X)-b(X_\tau) \> {\rm d } s \nonumber\\
&\quad+2p\ee\int_{0}^{t_n}|X-X_\tau|^{2p-2}
\< X-X_\tau,b(X_\tau)-b(X_\tau( \kappa_n)) \> {\rm d } s \nonumber\\
&\quad+	2p\ee\int_{0}^{t_n}|X-X_\tau|^{2p-2}
\< X-X_\tau,b(X_\tau( \kappa_n))-b_\tau(X_\tau( \kappa_n)) \> {\rm d } s. 
\end{align}
Then by the conditions \eqref{mon}, \eqref{b-btau}, the representation \eqref{b-taylor}, and the estimates \eqref{bound-xtau} and \eqref{err-R},  we obtain
 \begin{align*}
&\ee|X({t_n})-X_\tau({t_n})|^{2p} \\
&\le C\ee\int_{0}^{t_n}|X-X_\tau|^{2p}  {\rm d } s +C\ee\int_{0}^{t_n
}|b(X_\tau( \kappa_n))-b_\tau(X_\tau( \kappa_n))|^{2p}  {\rm d } s\\
&\quad+2p\ee\int_{0}^{t_n}|X-X_\tau|^{2p-2}
\< X-X_\tau, (s-\kappa_n) 
\nabla b(X_\tau(\kappa_n)) b_\tau(X_\tau(\kappa_n))\> {\rm d } s\\
&\quad+2p\ee\int_{0}^{t_n}|X-X_\tau|^{2p-2}\< X-X_\tau, 
\nabla b(X_\tau(\kappa_n)) {\rm B}(W-W(\kappa_n))\> {\rm d } s\\
&\quad+2p\ee\int_{0}^{t_n}|X-X_\tau|^{2p-2}
\< X-X_\tau, R(X_\tau(s), 
X_\tau(\kappa_n))\> {\rm d } s\\
&\le C\ee\int_{0}^{t_n}|X-X_\tau|^{2p} {\rm d } s+
\ee\int_{0}^{t_n}|R(X_\tau,X_\tau(\kappa_n))|^{2p} {\rm d } s \\
&\quad+ C \tau^{2p} \ee \int_{0}^{t_n} [(1+|X_\tau( \kappa_n)|^r)^p (L_4|X_\tau( \kappa_n)|^{q+1}+L_3)^{2p} \\
& \qquad  \qquad  \qquad + (1+|X_\tau(\kappa_n)|^q)^{2p}
(L_4|X_\tau(\kappa_n)|^{q+1}+L_3)^{2p}] {\rm d } s \\
&\quad+2p\ee\int_{0}^{t_n}|X-X_\tau|^{2p-2}
\< X-X_\tau,\nabla b(X_\tau(\kappa_n)) {\rm B}(W-W(\kappa_n))\> {\rm d } s
\\
&\le C\ee\int_{0}^{t_n}|X-X_\tau|^{2p}  {\rm d } s+e^{Ct_n}
\tau^{2p}(1+\ee|X_0|^{2p(2q+1)}+\ee|X_0|^{p(2q+2+r)}) + 2p \cdot I, 
\end{align*}
where 
\begin{align}\label{form-I}
I:= \ee\int_{0}^{t_n}|X-X_\tau|^{2p-2}
\< X-X_\tau, \nabla b(X_\tau(\kappa_n)){\rm B}(W-W(\kappa_n))\> {\rm d} s.
\end{align}

For the term $I$, we combine It\^o formula 
applied to the functional $|X-X_\tau|^{2p-2}$ with the integral representation of the SODE to get
\begin{align*}
 I&:=\ee \int_{0}^{t_n} |X(\kappa_n)-X_\tau(\kappa_n)|^{2p-2}\< X-X_\tau, \nabla b(X_\tau(\kappa_n)){\rm B}(W-W(\kappa_n))\> \dd s \\
&\quad  + 2(p-1) \ee \int_{0}^{t_n} \Big[ \int_{\kappa_n(s)}^s |X(r)-X_\tau(r)|^{2(p-2)} \< X(r)-X_\tau(r),  \\
&\qquad b(X(r))-b_\tau(X_\tau( \kappa_n(r))) \> {\rm d } r \Big]
\< X-X_\tau, \nabla b(X_\tau(\kappa_n)){\rm B}(W-W(\kappa_n))\> \dd s \\
&=\ee \int_{0}^{t_n} |X(\kappa_n)-X_\tau(\kappa_n)|^{2p-2}\< X(\kappa_n)-X_\tau(\kappa_n), \nabla b(X_\tau(\kappa_n)){\rm B}(W-W(\kappa_n))\> \dd s \\
&\quad + \ee \int_{0}^{t_n} |X(\kappa_n)-X_\tau(\kappa_n)|^{2p-2}\Big\< \int_{\kappa_n(s)}^s b(X(r))-b_\tau(X_\tau( \kappa_n(r)))\dd r ,\\
&\qquad
 \nabla b(X_\tau(\kappa_n))  {\rm B}(W-W(\kappa_n))\Big\> \dd s \\
&\quad  + 2(p-1) \ee \int_{0}^{t_n} \Big[ \int_{\kappa_n(s)}^s |X(r)-X_\tau(r)|^{2(p-2)} \< X(r)-X_\tau(r),  \\
&\qquad b(X(r))-b_\tau(X_\tau( \kappa_n(r))) \> {\rm d } r \Big]
\< X-X_\tau, \nabla b(X_\tau(\kappa_n)){\rm B}(W-W(\kappa_n))\> \dd s \\
&=:I_1+I_2+I_3.
\end{align*}
The term $I_1$ vanishes  because the Brownian increment  $W(s)-W(\kappa_n(s))$  is independent of  $\FFF_{\kappa_n(s)}$ for all $s \in (0, t)$.  For the remaining two terms,  using the H\"older and Young inequalities, we have
\begin{align*}
I_2+I_3 &\le C\int_{0}^{t_n}\int_{\kappa_n(s)}^s (\tau^{-1}\sup_{u\le s}\ee|X(u)-X_\tau(u)|^{2p})
 ^{\frac{p-1}{p}} \\
 &\quad \times [\tau^{p-1}\ee|b(X(r))-b_\tau(X(\kappa_n(r)))|^p
 |\nabla b(X_\tau(\kappa_n(s))) \\
 & \qquad \cdot {\rm B}(W(s)-W(\kappa_n(s)))|^p]^{\frac1p}
 {\rm d} r{\rm d} s \\
 &\le C\int_{0}^{t_n}\sup_{u\le s}\ee|X(u)-X_\tau(u)|^{2p}{\rm d} s  \\
& \quad  + C\tau^{p-1}\int_{0}^{t_n}\int_{\kappa_n(s)}^s 
 |b(X(r))-b_\tau(X(\kappa_n(r)))|^p \\
 &\qquad \times |\nabla b(X_\tau(\kappa_n(s))){\rm B}(W(s)-W(\kappa_n(s)))|^p{\rm d} r{\rm d} s \\
 &\le C\int_{0}^{t_n}\sup_{u\le s}\ee|X(u)-X_\tau(u)|^{2p}{\rm d} s \\
& \quad  + C\tau^{p-1}\int_{0}^{t_n}\int_{\kappa_n(s)}^s 
 (\ee|b(X(r))-b_\tau(X(\kappa_n(r)))|^{2p})^{\frac12} \\
 &\qquad \times(\ee|\nabla b(X_\tau(\kappa_n(s))){\rm B}
 (W(s)-W(\kappa_n(s)))|^{2p})^{\frac12}{\rm d} r{\rm d} s.
\end{align*}

Using \eqref{b-lip}, \eqref{b-btau}, \eqref{bound-x},
and \eqref{hol-x}, we have 
\begin{align*}
& \ee|b(X)-b_\tau(X(\kappa_n))|^{2p} \\
 &\le C\ee |b(X)-b(X(\kappa_n))|^{2p}
 +C\ee |b(X(\kappa_n))-b_\tau(X(\kappa_n))|^{2p} \nonumber\\
&\le C\ee\big[(1+|X|+|X(\kappa_n)|)^{2pq}|X-X(\kappa_n)|^{2p}\big]
\nonumber\\
&\quad+C \tau^{2p}\ee [(1+|X(\kappa_n)|^r)^p(L_4|X_\tau( \kappa_n)|^{q+1}+L_3)^{2p}] \nonumber\\
&\le C (\ee (1+ |X|+ |X(\kappa_n)|)^{p(3q+2)})^{\frac{2q}{3q+2}}(\ee |X-X(\kappa_n)|^{\frac{2p(3q+2)}{q+2}})^{\frac{q+2}{3q+2}} 
\nonumber\\
&\quad+C \tau^{2p} (1+ \ee |X(\kappa_n)|^{p(2q+2+r)}) \nonumber\\
&\le C e^{Ct_n} \tau^p(1+C t_n+\ee|X_0|^{p(3q+2)}+\ee|X_0|^{p(2q+2+r)}).
\end{align*}
Similarly,  using \eqref{b-add-poly} and \eqref{bound-xtau}, we obtain
\begin{align*}
\ee|\nabla b(X_\tau(\kappa_n)){\rm B}(W-W(\kappa_n))|^{2p}
&\le C\tau^p (1+\ee|X_\tau(\kappa_n)|^{2pq}) \nonumber\\
&\le C e^{Ct_n} \tau^p (1+C t_n+\ee|X_0|^{2pq}).
\end{align*}
In combination with the above two estimates, we have
\begin{align*}
&\ee|X({t_n})-X_\tau({t_n})|^{2p} \\
&\le C\int_{0}^{t_n} \sup_{u\le s}
|X-X_\tau|^{2p}{\rm d }s+
C e^{Ct_n}\tau^{2p}(1+C t_n+\ee|X_0|^{p(2q+2+r)}+\ee|X_0|^{2p(2q+1)}),
\end{align*}
from which we conclude the desired result \eqref{err-L2*} by Gronwall inequality.
 \end{proof}

\begin{remark}
The strong order one convergence result \eqref{err-L2*}, in combination with the arguments developed in \cite{NNZ25}, yields the time-independent weak order one convergence under Wasserstein distance for the TEM scheme \eqref{tem-} towards Eq. \eqref{sde}, under the following additional convex at infinity condition, consistent with our main Assumptions \ref{ap}, \ref{ap-non}, and \ref{ap+}, with some positive constants ${\widehat L}_1$, ${\widehat L}_2$, and ${\widehat q} \in [0, q)$:
 		\begin{align*} 
\<x-y,b(x)-b(y)\> & \leq -{\widehat L}_1 |x-y|^2(|x|^q+|y|^q)+{\widehat L}_2 |x-y|^2(|x|^{\widehat q}+|y|^{\widehat q}), \quad  x,y \in \mathbb{ R }^d.
 		\end{align*} 
Then the time-independent convergence rate between the numerical invariant measure of the TEM scheme \eqref{tem-} and the exact invariant measure of Eq. \eqref{sde} follows immediately.
 The details for this and the multiplicative noise case have been developed in \cite{LWWZ25}.  
		
\end{remark}

\section{Numerical Experiments}
\label{sec4}

  In this section, we present two numerical experiments to verify our theoretical results: Theorem \ref{tm-erg} in Section \ref{sec2} and Theorems \ref{tm-err} and \ref{tm-err+} in Section \ref{sec3}, respectively.

 The first numerical test is given to Eq. \eqref{sde} with $d=m=1$, $b(x)=(1-x^2) x$, and $\sigma(x)=(1+x^2)/2$, $x \in \rr$, i.e., \eqref{bs-ex} with $k=1$, $a_3=-1, a_2=0, a_1=1, a_0=0$, and $c_2=c_0=1/2, c_1=0$:
 \begin{align}\label{ac}  
    \dd X(t) = (1-X(t)^2) X(t) \dd t + 0.5 (1+X(t)^2) \dd W(t), ~ t \ge 0; 
    \quad X(0)=X_0,  
\end{align}   
so that the TEM scheme \eqref{tem} with \eqref{bs-tau} becomes 
\begin{align} \label{tem-ex1}
 	Z_n=Z_{n-1}+
 \frac{(1-Z_{n-1}^2) Z_{n-1} \tau}{(1+\tau|Z_{n-1}|^4)^{1/2}}
	+\frac{(1+Z_{n-1}^2) \delta_{n-1} W}{2 (1+\sqrt{\tau}|Z_{n-1}|^2)^{1/2}},
~ n \in \nn_+; \quad Z_0=X_0.
 \end{align}  
In this case, Assumption \ref{ap-non} is trivially satisfied, and the condition \eqref{bs-ex-con} holds with $2a_3+c_2^2=-7/4<0$.
Consequently, Assumptions \ref{ap} (with $q=2$ and $L_2 \in (0, 7/4)$) and \ref{ap+} (with $2 = p<p*<3.5$) follow from Remarks \ref{rk-ap}(2) and \ref{rk-ap+}(1). 
Consequently, the TEM scheme \eqref{tem-ex1} is geometrically ergodic for any small stepsize $\tau$ (e.g., $\tau<49/128$) and converges to Eq. \eqref{ac} with strong order $1/2$ (under the $L^4(\Omega)$-norm), according to Theorem \ref{tm-erg}(1) and  Theorem \ref{tm-err}, respectively. 

\begin{figure}[h]
    \centering 
    \includegraphics[width=1\textwidth]{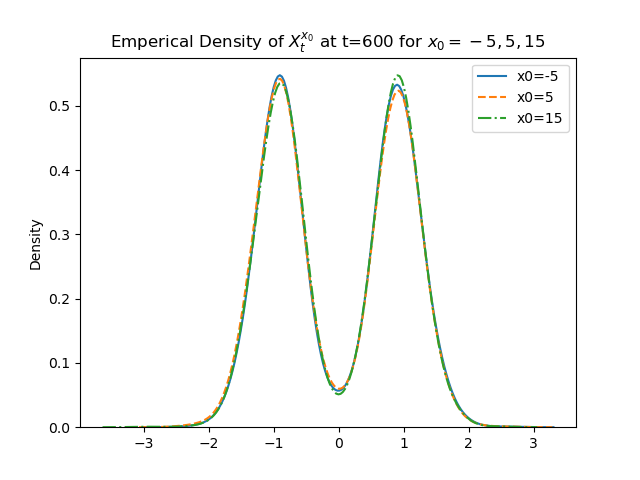}
    \caption{Empirical density of \eqref{tem} for \eqref{sde}}\label{fig-mul}
  \end{figure}

\begin{figure}[h]
    \centering 
    \includegraphics[width=1\textwidth]{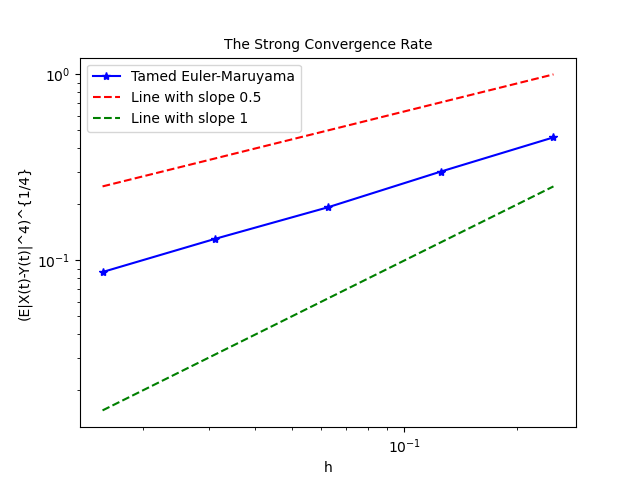}
    \caption{Strong error estimates of \eqref{tem} for \eqref{sde}}\label{fig-mul+}
  \end{figure}

We take $\tau=0.3<49/128$ and initial data $X_0=-5, 5, 15$ for the numerical experiment.
 It is clear from Figure \ref{fig-mul} that the empirical density functions with different initial data obtained via kernel density estimation at $n=2000$ (corresponding to $t=600$) become nearly indistinguishable. This convergence to a common distribution indicates the mixing property (yielding ergodicity) of \eqref{tem-ex1} and thus corroborates the theoretical result in Theorem \ref{tm-erg}.  
To demonstrate the strong convergence order, we set $ T=16$ and $ X_0=1$, and identify the reference exact solution with a numerical solution generated by \eqref{tem-ex1} using a fine time stepsize $\tau=2^{-8}$. 
The other numerical approximations are calculated by \eqref{tem-ex1} with
five different time stepsizes  $\tau = 2^{-i},i=2,3,4,5,6$.
Then $M=5000$ sample paths will be used to simulate the expectation.
Figure \ref{fig-mul+} shows that the TEM scheme \eqref{tem-ex1} converges to Eq. \eqref{ac} with strong order $1/2$ (in the $L^4(\Omega)$-norm), thereby verifying the theoretical result in Theorem \ref{tm-err}.

  We also simulate the additive noise case by considering Eq. \eqref{sde} with $b(x)=(1-x^2) x$ and $\sigma(x)=1/2$, $x \in \rr$:
 \begin{align}\label{ac-}  
\dd X(t) = (1-X(t)^2) X(t) \dd t + 0.5 \dd W(t), ~ t \ge 0; 
    \quad X(0)=X_0,  
\end{align}   
so that the TEM scheme \eqref{tem-} with \eqref{bs-tau} becomes 
\begin{align} \label{tem-ex2}
Z_n=Z_{n-1}+
 \frac{(1-Z_{n-1}^2) Z_{n-1} \tau}{(1+\tau|Z_{n-1}|^4)^{1/2}}
	+0.5\delta_{n-1} W,
~ n \in \nn_+; \quad Z_0=X_0.
 \end{align}   
 In this case, Assumption \ref{ap-non} is clearly satisfied, and conditions  \eqref{b-grow}, \eqref{b-coe}, \eqref{s-grow}, and \eqref{b-add-poly} hold (with $q=2$ and $L_2' \in (0, 2)$).
By Theorem \ref{tm-erg}(2) and Theorem \ref{tm-err+}, the TEM scheme \eqref{tem-ex2} is geometrically ergodic for any small stepsize $\tau<1/2$ and converges to  Eq. \eqref{ac-} with strong order $1$ (under the $L^4(\Omega)$-norm). 

\begin{figure}[h]
    \centering 
    \includegraphics[width=1\textwidth]{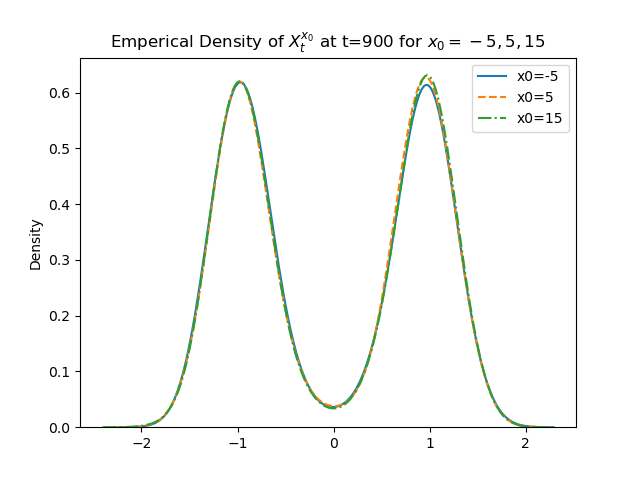}
    \caption{Empirical density of \eqref{tem*} for \eqref{sde}}\label{fig-add}
  \end{figure}

\begin{figure}[h]
    \centering 
    \includegraphics[width=1\textwidth]{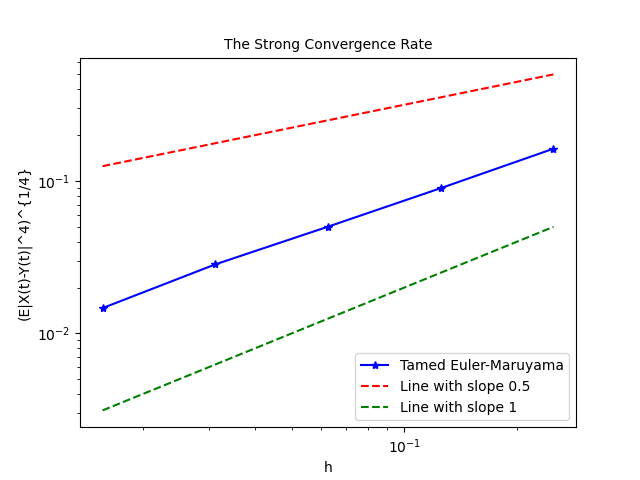}
    \caption{Strong error estimates of \eqref{tem*} for \eqref{sde}}\label{fig-add+}
  \end{figure}

Similarly to the first test, we use the stepsize $\tau=0.45$ and the initial data $ X_0 = -5, 5, 15$. 
It is clear from Figure \ref{fig-add} that the empirical density functions with different initial data plotted by kernel density estimation at $n=2000$ (corresponding to $t=900$) become nearly indistinguishable, which indicates the mixing property of \eqref{tem-ex2} and thus verifies the theoretical result in Theorem \ref{tm-erg}.   
To demonstrate the strong convergence order, we set $ T=16$ and $ X_0=1$, and use step sizes and the number of sample paths similar to those in the first test.
Figure \ref{fig-add+} shows that the TEM scheme \eqref{tem-ex2} converges to  Eq. \eqref{ac-}  with strong order $1$ (in the $L^4(\Omega)$-norm) and verifies the theoretical result in Theorem \ref{tm-err+}.

 \section*{Declarations}

{\bf Conflict of interest}. The authors have no competing interests to declare relevant to this article's content.
The data are available from the corresponding author on reasonable request.

 \section*{Acknowledgements}

We thank the anonymous referees for their helpful comments and suggestions. 
The authors are supported by the National Natural Science Foundation of China, No. 12101296, Guangdong Basic and Applied Basic Research Foundation, No. 2024A1515012348, and Shenzhen Basic Research Special Project (Natural Science Foundation) Basic Research (General Project), Nos. JCYJ20220530112814033 and JCYJ20240813094919026.
    
  \bibliographystyle{amsplain}
  \bibliography{bib.bib}

\end{document}